\documentclass[12pt]{article}
\usepackage{amssymb,amsmath}
\usepackage{cases}
\usepackage{amsfonts}
\usepackage{color,xcolor}
\usepackage[left=2.0cm,right=2.0cm,top=2.0cm,bottom=2.0cm]{geometry}
\usepackage[colorlinks,citecolor=blue,urlcolor=blue]{hyperref}

\newtheorem{theorem}{Theorem}[section]
\newtheorem{corollary}{Corollary}[section]
\newtheorem{lemma}{Lemma}[section]

\newtheorem{remark}{Remark}[section]

\newcommand{\bal}{\begin{align}}
\newcommand{\bbal}{\begin{align*}}
\newcommand{\beq}{\begin{equation}}
\newcommand{\eeq}{\end{equation}}
\newcommand{\bca}{\begin{cases}}
\newcommand{\eca}{\end{cases}}
\def\div{\mathord{{\rm div}}}
\newcommand{\pa}{\partial}
\newcommand{\fr}{\frac}
\newcommand{\na}{\nabla}

\newcommand{\cd}{\cdot}
\newcommand{\ep}{\varepsilon}
\newcommand{\dd}{\mathrm{d}}

\newcommand{\R}{\mathbb{R}}

\newcommand{\les}{\lesssim}

\newcommand{\bi}{\Big}

\begin{document}
\title{Global smooth solutions of the damped Boussinesq equations with a class of large initial data}

\author{Jinlu Li$^{1}$\footnote{E-mail: lijinlu@gnnu.cn}, Xing Wu$^{2}$\footnote{E-mail: ny2008wx@163.com}, Weipeng Zhu$^{3}$\footnote{E-mail:  mathzwp2010@163.com}\\
\small $^1$\it School of Mathematics and Computer Sciences, Gannan Normal University, \\
\small Ganzhou, Jiangxi, 341000, China\\
\small $^2$\it College of Information and Management Science,
Henan Agricultural University,\\
\small Zhengzhou, Henan, 450002, China\\
\small $^3$\it  School of Mathematics and Information Science, Guangzhou University, Guangzhou 510006, China}

\date{}

\maketitle\noindent{\hrulefill}

{\bf Abstract:} The global regularity problem concerning the  inviscid Boussinesq equations remains an open problem. In an attempt to understand this problem, we examine the damped  Boussinesq equations and study
how damping affects the regularity of solutions. In this paper, we consider the global existence to the damped Boussinesq equations with a class of large initial data, whose $B^{s}_{p,r}$ or $\dot{B}^{s}_{p,r}$ norms can be arbitrarily large. The idea is splitting the linear Boussinesq equations from the damped Boussinesq equations, the exponentially decaying solution of the former equations together with the structure of the Boussinesq equations help us to obtain the global smooth solutions.


\vskip0mm\noindent{\hrulefill}

\section{Introduction}\label{sec1}
This paper considers the global smooth solutions for the incompressible  Boussinesq equations with damping
\begin{eqnarray}\label{boussineq}
        \left\{\begin{array}{ll}
          \partial_tu+u\cd\na u+\nu u+\na p=\theta e_d,& x\in \R^d,t>0,\\
          \partial_t\theta+u\cd\na \theta +\lambda \theta=0,& x\in \R^d,t>0,\\
         \div u=0,& x\in \R^d,t\geq0,\\
          (u,\theta)|_{t=0}=(u_0,\theta_0),& x\in \R^d,\end{array}\right.
        \end{eqnarray}
where $u$ is the velocity vector, $p$ is the pressure, $\theta$ denotes the temperature or density which is a scalar function. $e_d=(0, 0,...,0,1)^T$ and $\nu, \lambda$ are two positive parameters, standing for kinematic viscosity and thermal diffusivity, respectively.

The Boussinesq equations model large scale atmospheric and oceanic flows that are
responsible for cold fronts and the jet stream \cite{Majda 2003, Pedlosky 1987} and  mathematically has received significant attention, since it has a vortex stretching effect similar to that in the 3D incompressible flow. When $\nu u$ is replaced by $-\nu \Delta u$, $\lambda \theta$  by $-\lambda \Delta \theta$, \eqref{boussineq} becomes the standard viscous Boussinesq equations, the global in time regularity in two dimension is well understood even in the zero diffusivity ($\nu >0, \lambda=0$) or the zero viscosity case($\nu =0, \lambda>0$)  \cite{Abidi 2007, Chae 2006, Chae 1997, Chen 2017, Hmidi 2009, Hou 2005, Lai 2011}, however the global regularity in dimension
three appears to be out of reach. while $\nu=0$ and $\lambda=0$, (\ref{boussineq}) is reduced to the inviscid  Boussinesq equations, due to the absence of dissipative terms, the global solution or finite-time singularity  evoluting from general initial data remains unsolved in spite of the progress on the local well-posedness and regularity criteria \cite{Cui 2012, Hassainia 2015, Liu 2010, Taniuchi 2002, Tao 2017, Wan 2016}.
Recently, following the convex integration method, Tao and Zhang \cite{Tao 2017} obtained the  H\"{o}lder continuous solution with compact support both in space and time for  inviscid 2D Boussinesq equations.

When adding velocity damping term $\nu u$ and temperature damping term $\lambda \theta$, Adhikar et al. \cite{Adhikar 2014} proved (\ref{boussineq}) admits  a unique global small solution  with the initial data satisfying
$$\|\nabla u_0\|_{\dot{B}_{\infty, 1}^0}<\min\{\frac{\nu}{2C}, \frac{\lambda}{C}\}, \qquad \|\nabla \theta_0\|_{\dot{B}_{\infty, 1}^0}<\frac{\nu}{2C}\|\nabla u_0\|_{\dot{B}_{\infty, 1}^0},$$
where $\dot{B}_{\infty, 1}^0$  is the homogeneous Besov space. We are interested in that  whether or nor (\ref{boussineq}) possesses a global  solution without the  smallness assumption. If not, it may be helpful to obtain the global solutions for a class of large initial data. Our effort in this paper is precisely based on this. In what follows, let $(U,\Theta)$ be the solution of the "linearized damped Boussinesq equation",
\begin{eqnarray}\label{2}
        \left\{\begin{array}{ll}
          \pa_t U+\nu U+\na p=\Theta e_d,\\
          \pa_t \Theta+\lambda\Theta=0,\\
          \div U=0,\\
          (U,\Theta)|_{t=0}=(U_0,\Theta_{0}).\end{array}\right.
\end{eqnarray}
The second equation in (\ref{2}) satisfied by $ \Theta$ is a simple linear equation, and has the solution $\Theta=e^{-\lambda t}\Theta_0$. If we resort to the equations (\ref{2}) of the vorticity $W=\nabla\times U$, then $W$ satisfys
\begin{eqnarray*}
&&2D:  \qquad \partial_tW+\nu W=\partial_1\Theta,\\
&&3D:  \qquad \partial_tW+\nu W=(\partial_2\Theta, -\partial_1\Theta, 0)^T.
\end{eqnarray*}
Therefor both $W$ and $ \Theta$  have explicit exponential decay expressions, so does for $U$. By virtue of the good exponential decay property of $U$ and $ \Theta$, we obtain the global solution of (\ref{boussineq}) which is constructed with the large initial data.

We first recall the definition of the Besov spaces. Choose a radial, non-negative, smooth and radially decreasing function $\chi: {\mathbb R}^d\to [0, 1]$ such that it is supported in $\{\xi\in \R^d:|\xi|\leq \frac43\}$ and $\chi\equiv 1$ for $|\xi|\leq \frac34$. Let $\varphi(\xi)=\chi(\frac{\xi}{2})-\chi(\xi)$. Then $\varphi$ is supported in the ring $\{\xi\in\mathbb{R}^d:\frac 34\leq|\xi|\leq \frac83\}$ and $\varphi\equiv 1$ for $\frac43\leq |\xi|\leq \frac32$. For $u \in \mathcal{S}'$, $q\in {\mathbb Z}$, we define the Littlewood-Paley operators: $\dot{\Delta}_q{u}=\mathcal{F}^{-1}(\varphi(2^{-q}\cdot)\mathcal{F}u)$, ${\Delta}_q{u}=\dot{\Delta}_q{u}$ for $q\geq 0$, ${\Delta}_q{u}=0$ for $q\leq -2$ and $\Delta_{-1}u=\mathcal{F}^{-1}(\chi \mathcal{F}u)$, and $S_q{u}=\mathcal{F}^{-1}\big(\chi(2^{-q}\xi)\mathcal{F}u\big)$.
Here we use ${\mathcal{F}}(f)$ or $\widehat{f}$ to denote
the Fourier transform of $f$.

The standard vector-valued functions $u:{\mathbb R}^d\to {\mathbb R}^d$ in Besov spaces $B^s_{p,r}$ and $\dot B^s_{p,r}$ can be defined by
\begin{align*}
\|u\|_{B^s_{p,r}}&\triangleq \big|\big|(2^{js}\|\Delta_j{u}\|_{L^p})_{j\in {\mathbb Z}}\big|\big|_{\ell^r}<+\infty,\\
\|u\|_{\dot{B}^s_{p,r}}&\triangleq \big|\big|(2^{js}\|\dot{\Delta}_j{u}\|_{L^p})_{j\in {\mathbb Z}}\big|\big|_{\ell^r}<+\infty.
\end{align*}
Obviously, if $\mathrm{supp} \ \hat{u}\in \{\xi:\frac43\leq |\xi|\leq \frac32\}$, then we have
\bbal
||u||_{L^p}=||u||_{B^s_{p,r}}=||u||_{\dot{B}^s_{p,r}}.
\end{align*}

Our main result is stated as follows.
\begin{theorem}\label{the1.1} Let $d=2, 3$. Assume that the initial data fulfills ${\rm{div}}u_0=0$ and
$$u_0=U_0+v_0\quad \mbox{and}\quad \theta_0=\Theta_0+\vartheta_0$$
with
\begin{eqnarray}\label{Equ1.2}
\mathrm{supp} \ \hat{U}_0(\xi), \ \hat{\Theta}_0(\xi) \subset\mathcal{C}:=\Big\{\xi \big| \  \frac43\leq |\xi|\leq \frac32\Big\} .
\end{eqnarray}
Denote
\bbal
E_0=\int^{\infty}_{0}(||U\cd\na U||_{H^3}+||U\cd\na\Theta||_{H^3})\dd t, \quad F_0=\int^{\infty}_{0}||(U, \Theta)||_{L^\infty}\dd t.
\end{align*}
There exists a sufficiently small positive constant $\delta$, and a universal constant $C$ such that if
\begin{align}\label{condition}
\Big(||v_0||^2_{H^3}+||\vartheta_0||^2_{H^3}+E_0\Big)\cd\exp\Big( CF_0+CE_0\Big)\leq \delta,
\end{align}
then the system \eqref{boussineq} has a unique global solution.
\end{theorem}

\begin{corollary}\label{cor1.1}
For $d=2$, we let $v_0=\vartheta_0=0$ and
\begin{eqnarray*}
&U_0=\na^{\bot}a_0=
\begin{pmatrix}
\pa_2a_0 \\ -\pa_1a_0
\end{pmatrix}
,\qquad \Theta_0=a_0,
\end{eqnarray*}
where
\bbal
\mathrm{supp} \ \hat{a}_0(\xi)\subset\mathcal{\widetilde{C}}:=\Big\{\xi \big| \ |\xi_1-\xi_2|\leq \ep, \ \frac43\leq |\xi|\leq \frac32\Big\}.
\end{align*}
Then, direct calculations show that the left side of (\ref{condition}) becomes
\bbal
 C\ep||a_0||_{L^2}||\hat{a}_0||_{L^1}\mathrm{exp}
 \Big({C\ep||a_0||_{L^2}||\hat{a}_0||_{L^1}+C||\hat{a}_0||_{L^1}}\Big).
\end{align*}
For $d=3$, we let $v_0=\vartheta_0=0$ and
\begin{eqnarray*}
&U_0=
\begin{pmatrix}
\pa_2a_0 \\ -\pa_1a_0, \\ 0
\end{pmatrix}
,\qquad \Theta_0=a_0,
\end{eqnarray*}
where
\bbal
\mathrm{supp} \ \hat{a}_0(\xi)\subset\mathcal{\widetilde{C}}_0:=\Big\{\xi \big| \ |\xi_1-\xi_2|\leq \ep, \ \frac{41}{30}\leq |\xi_h|\leq \frac{22}{15}, \ \ep\leq |\xi_3|\leq 2\ep\Big\}.
\end{align*}
Then, direct calculations show that the left side of (\ref{condition}) becomes
\bbal
 C\ep||a_0||_{L^2}||\hat{a}_0||_{L^1}\mathrm{exp}\Big({C\ep||a_0||_{L^2}||\hat{a}_0||_{L^1}+C||\hat{a}_0||_{L^1}}\Big).
\end{align*}
\end{corollary}

\begin{remark}
For $d=2$, we set
\bbal
a_0(x_1,x_2)=\ep^{-1}\bi(\log\log\frac1\ep\bi)^{\frac12}  \chi(x_1,x_2),
\end{align*}
where the smooth functions $\chi$ satisfying $\hat{\chi}(-\xi_1,-\xi_2)=\hat{\chi}(\xi_1,\xi_2)$,
\begin{align*}
\mathrm{supp} \hat{\chi}\subset \mathcal{\widetilde{C}},\quad \hat{\chi}(\xi)\in[0,1]\quad\mbox{and} \quad \hat{\chi}(\xi)=1 \quad\mbox{for} \quad \xi\in\mathcal{\widetilde{C}}_1,
\end{align*}
where
\begin{align*}
\mathcal{\widetilde{C}}_1\triangleq\Big\{\xi\in\R^2: \ |\xi_1-\xi_2|\leq \frac12\ep,\ \frac{25}{18}\leq |\xi|\leq \frac{13}{9}\Big\}.
\end{align*}
Then, direct calculations show that the left side of \eqref{condition} becomes
\begin{align*}
C\ep^{\frac12}\bi(\log\log \frac1\ep\bi)\exp\bi(C\log\log \frac1\ep\bi).
\end{align*}
In fact, one has
\begin{align*}
||\hat{a}_0||_{L^1}\approx \Big(\log\log\frac1\ep\Big)^\frac12\quad\mbox{and}\quad||{a}_0||_{L^2}\approx \ep^{-\fr12}\Big(\log\log\frac1\ep\Big)^\frac12.
\end{align*}
Therefore, choosing $\ep$ small enough, we deduce that the system (\ref{boussineq}) has a unique global solution. Moreover, we can show that
$$||u_0||_{L^\infty}\gtrsim \bi(\log\log \frac1\ep\bi)^\frac12,\qquad ||\theta_0||_{L^\infty}\gtrsim \bi(\log\log \frac1\ep\bi)^\frac12.$$
\end{remark}

\begin{remark}\label{rem1.2}
For $d=3$ and $1<p< \infty$ we set
$$a_0(x_1,x_2,x_3)=\ep^{-\frac{2(p-1)}{p}}(\log\log\frac1\ep)^{\frac12}  \chi(x_1,x_2)\phi(x_3).$$
where the smooth functions $\chi,\phi$ satisfying
\begin{align*}
\mathrm{supp} \hat{\chi}\in \mathcal{\widetilde{C}},\quad \hat{\chi}(\xi)\in[0,1]\quad\mbox{and} \quad \hat{\chi}(\xi)=1 \quad\mbox{for} \quad \xi\in\mathcal{\widetilde{C}}_1,
\end{align*}
and
\begin{align*}
\mathrm{supp} \hat{\phi}(\xi') \in [\ep,2\ep],\quad \hat{\phi}(\xi') \in [0,1] \quad\mbox{and} \quad \hat{\phi}(\xi)=1 \quad\mbox{for}  \quad \xi'\in[\frac54\ep,\frac74\ep].
\end{align*}
Then, direct calculations show that the left side of \eqref{condition} can be bounded by
\begin{align*}
C\ep^{\frac{4}{p}} \log\log \frac1\ep e^{C\log\log\frac1\ep}.
\end{align*}
Therefore, choosing $\ep$ small enough, we deduce that the system (\ref{boussineq}) has a global solution. $$||u_0||_{L^p}\gtrsim \bi(\log\log \frac1\ep\bi)^\frac12,\qquad ||\theta_0||_{L^p}\gtrsim \bi(\log\log \frac1\ep\bi)^\frac12.$$
\end{remark}

{\bf Notations}: Let $\beta=(\beta_1,\beta_2,\beta_3)\in \mathbb{N}^3$ be a multi-index and $D^{\beta}=\pa^{|\beta|}/\pa^{\beta_1}_{x_1}\pa^{\beta_2}_{x_2}\pa^{\beta_3}_{x_3}$ with $|\beta|=\beta_1+\beta_2+\beta_3$. For the sake of simplicity, $a\lesssim b$ means that $a\leq Cb$ for some ``harmless" positive constant $C$ which may vary from line to line. $[A,B]$ stands for the commutator operator $AB-BA$, where $A$ and $B$ are any pair of operators on some Banach space $X$. We also use the notation $||f_1,\cdots,f_n||_{X}\triangleq||f_1||_{X}+\cdots+||f_n||_{X}$.

\section{Proof of Theorem \ref{the1.1}}

In this section, we give the proof of Theorem \ref{the1.1}. Before giving the proof, we present some estimates which will be used in the proof of Theorem \ref{the1.1}.

{\bf Proof of Theorem \ref{the1.1}}\quad  Denoting $v=u-U$ and $\vartheta=\theta-\Theta$, we can reformulate the system \eqref{boussineq} and \eqref{2} equivalently as
\begin{eqnarray}\label{5}
        \left\{\begin{array}{ll}
\partial_tv+v\cd\na v+U\cd\na v+v\cd\na U+\nu v+\na p'-\vartheta e_d=-U\cd\na U:=f,\\
\partial_t\vartheta+v\cd\na \vartheta+U\cd\na\vartheta+v\cd\na\Theta+\lambda \vartheta=-U\cd\na\Theta:=g,\\
\div v=0,\\
(v,c)|_{t=0}=(v_0,c_0).\end{array}\right.
\end{eqnarray}
Applying $D^\beta$ to $\eqref{5}_1$ and $\eqref{5}_2$ respectively and taking the scalar product  with $\sigma D^\beta v$ and $D^\beta \vartheta$ respectively, adding them together and then summing the result over $|\beta|\leq 3$, we get
\bal\label{z0}
\fr12\frac{\dd}{\dd t}\Big(\sigma ||v||^2_{H^3}+||\vartheta||^2_{H^3}\Big)+\sigma \nu|| v||^2_{H^3}+\lambda||\vartheta||^2_{H^3}-\sigma(\vartheta, v_d)_{H^3}\triangleq\sum^{4}_{i=1}I_i,
\end{align}
where
\bbal
&I_1=-\sigma\sum_{0<|\beta|\leq 3}\int_{\R^3}[D^{\beta},v\cd] \na v\cd D^\beta v\dd x
-\sum_{0<|\beta|\leq 3}\int_{\R^3}[D^{\beta},v\cd] \na \vartheta\cd D^\beta \vartheta\dd x,
\\&I_2=-\sigma\sum_{0<|\beta|\leq 3}\int_{\R^3}D^{\beta}(U\cd \na v)\cd D^\beta v\dd x-\sum_{0<|\beta|\leq 3}\int_{\R^3}D^{\beta}(U\cd \na \vartheta)\cd D^\beta \vartheta\dd x,
\\&I_3=-\sigma\sum_{0\leq|\beta|\leq 3}\int_{\R^3}D^{\beta}(v\cd \na U)\cd D^{\beta}v\dd x-\sum_{0\leq|\beta|\leq 3}\int_{\R^3}D^{\beta}(v\cd \na \Theta)\cd D^{\beta}\vartheta\dd x,
\\&I_{4}=-\sigma\sum_{0\leq |\beta|\leq 3}\int_{\R^3}D^{\beta}(U\cd\na U)\cd D^{\beta}v\dd x-\sum_{0\leq |\beta|\leq 3}\int_{\R^3}D^{\beta}(U\cd\na\Theta)\cd D^{\beta}\vartheta\dd x.
\end{align*}
Next, we need to estimate the above terms one by one.

According to the commutate estimate (See \cite{Majda 2001}),
\bal\label{l0}
\sum_{|\alpha|\leq m}||[D^{\alpha},\mathbf{g}]\mathbf{f}||_{L^2}\leq C(||\mathbf{f}||_{H^{m-1}}||\na \mathbf{g}||_{L^\infty}+||\mathbf{f}||_{L^\infty}||\mathbf{g}||_{H^m}),
\end{align}
we obtain
\bal
I_1\leq&~\sigma\sum_{0<|\beta|\leq 3}||[D^{\beta},v\cd] \na v||_{L^2}||\na v||_{H^2}+\sum_{0<|\beta|\leq 3}||[D^{\beta},v\cd]\na \vartheta||_{L^2}||\na\vartheta||_{H^2}\nonumber\\
\leq&~C||\na v||_{L^\infty}||v||_{H^3}||\na v||_{H^2}+C(||\na v||_{L^\infty}||\na \vartheta||_{H^2}+||v||_{H^3}||\na \vartheta||_{L^\infty})||\na \vartheta||_{H^2}\nonumber\\
\leq&~C||v||_{H^3}\Big(||\na v||^2_{H^2}+||\na \vartheta||^2_{H^2}\Big)\leq C||v||_{H^3}\Big(||v||^2_{H^3}+||\vartheta||^2_{H^3}\Big).\label{z1}
\end{align}
Invoking the following calculus inequality which is just a consequence of Leibniz's formula,
\bbal
\sum_{|\beta|\leq 3}||[D^{\beta},\mathbf{g}]\mathbf{f}||_{L^2}\leq C(||\na \mathbf{g}||_{L^\infty}+||\na^3 \mathbf{g}||_{L^\infty})||\mathbf{f}||_{H^2},
\end{align*}
we obtain
\bal
I_2\leq&~\sigma\sum_{0<|\beta|\leq 3}||[D^{\beta},U\cd] \na v||_{L^2}||\na v||_{H^2}+\sum_{0<|\beta|\leq 3}||[D^{\beta},U\cd] \na \vartheta||_{L^2}||\na \vartheta||_{H^2}\nonumber\\
\leq&~C\Big(||\na U||_{L^\infty}+||\na^3 U||_{L^\infty}\Big)\Big(||v||^2_{H^3}+||\vartheta||^2_{H^3}\Big)\leq C||U||_{L^\infty}\Big(||v||^2_{H^3}+||\vartheta||^2_{H^3}\Big),\label{z4}
\end{align}
where we have used the fact that $\mathrm{supp} \ \hat{U}(\xi), \ \mathrm{supp} \ \hat{\Theta}(\xi) \subset \mathrm{supp} \ \hat{U}_0(\xi) \cup \mathrm{supp} \ \hat{\Theta}_0(\xi) \subset\mathcal{C}$.
By Leibniz's formula and H\"{o}lder's inequality, one has
\bal
I_3\leq&~\sigma||v\cd \na U||_{H^3}||v||_{H^3}+||v\cd \na \Theta||_{H^3}||\vartheta||_{H^3}\nonumber\\
\leq&~ C\Big(||\na (U, \Theta)||_{L^\infty}+||\na^4 (U, \Theta)||_{L^\infty}\Big)\Big(||v||^2_{H^3}+||\vartheta||^2_{H^3}\Big)\nonumber\\
 \leq&~ C||(U, \Theta)||_{L^\infty}\Big(||v||^2_{H^3}+||\vartheta||^2_{H^3}\Big).\label{z8}
\end{align}
Using  H\"{o}lder's inequality and Young inequality, we deduce
\bal
I_{4}&\leq~ C(||U\cd\na U||_{H^3}+||U\cd\na\Theta||_{H^3})(||v||_{H^3}+||\vartheta||_{H^3})\nonumber
\\&\leq~ C||(U\cd\na U, U\cd\na\Theta)||_{H^3}+C||(U\cd\na U, U\cd\na\Theta)||_{H^3}(||v||^2_{H^3}+||\vartheta||^2_{H^3}).\label{z10}
\end{align}
Putting all the estimates \eqref{z1}--\eqref{z10} into \eqref{z0}, we obtain
\bal\label{z12}
&\quad \frac{\dd}{\dd t}\Big(\sigma ||v||^2_{H^3}+||\vartheta||^2_{H^3}\Big)+\sigma \nu|| v||^2_{H^3}+\lambda||\vartheta||^2_{H^3}-\sigma (\vartheta,v_d)_{H^3}\nonumber\\&\les
||v||_{H^3}\Big(||v||^2_{H^3}+||\vartheta||^2_{H^3}\Big)
+\Big(||(U, \Theta)||_{L^\infty}+||(U\cd\na U, U\cd\na\Theta)||_{H^3}\Big)\Big(||v||^2_{H^3}+||\vartheta||^2_{H^3}\Big)\nonumber\\
&\quad+||(U\cd\na U, U\cd\na\Theta)||_{H^3}.
\end{align}
According to the inequality
\bbal
\sigma (\vartheta,v_d)_{H^3}\leq C\sigma ^{\frac32}||v||^2_{H^3}+C\sigma ^{\frac12}||\vartheta||^2_{H^3},
\end{align*}
then we can choose $\sigma $ small enough such that
\bbal
\sigma \nu|| v||^2_{H^3}+\lambda||\vartheta||^2_{H^3}-\sigma (\vartheta,v_d)_{H^3}
\approx \sigma \nu|| v||^2_{H^3}+\lambda||\vartheta||^2_{H^3}.
\end{align*}
For simplicity, we  denote
\bbal
A(t)=\sigma ||v||^2_{H^3}+||\vartheta||^2_{H^3}, \qquad B(t)=\sigma\nu|| v||^2_{H^3}+\lambda||\vartheta||^2_{H^3},
\end{align*}
then (\ref{z12}) can be rewritten as
\bal\label{z13}
\frac{\dd}{\dd t}A(t)+B(t)&\leq CA^{\frac12}(t)B(t)+C\Big(||(U, \Theta)||_{L^\infty}+||(U\cd\na U, U\cd\na\Theta)||_{H^3}\Big)A(t)\nonumber\\
&\quad+C||(U\cd\na U, U\cd\na\Theta)||_{H^3}.
\end{align}
Now, we define
\bbal
\Gamma\triangleq\sup\{t\in[0,T^*):\sup_{\tau\in[0,t]}A(\tau)\leq \eta\},
\end{align*}
where $\eta$ is a small enough positive constant which will be determined later. Assume that $\Gamma<T^*$. For all $t\in[0,\Gamma]$, we obtain from \eqref{z13} that
\bal
\frac{\dd}{\dd t}A(t)+B(t)&\leq C\Big(||(U, \Theta)||_{L^\infty}+||(U\cd\na U, U\cd\na\Theta)||_{H^3}\Big)A(t)+C||(U\cd\na U, U\cd\na\Theta)||_{H^3}.
\end{align}
which together with the assumption (\ref{condition}) yields that
\bbal
A(t)\leq &C\Big(||v_0||^2_{H^3}+||\vartheta_0||^2_{H^3}+E_0\Big)\cd\exp\Big( CE_0+CF_0\Big)\leq C\delta.
\end{align*}
Choosing $\eta=2C\delta$, thus we can get $\sup_{\tau\in[0,t]}A(\tau)\leq \fr\eta2$ for $t\leq \Gamma$. So if $\Gamma<T^*$, due to the continuity of the solutions, we can obtain there exists $0<\epsilon\ll1$ such that $\sup_{\tau\in[0,t]}A(\tau)\leq \eta$ for $t\leq \Gamma+\epsilon<T^*$, which contradicts with the definition of $\Gamma$. Thus, we can conclude $\Gamma=T^*$ and
\bbal
\sup_{\tau\in[0,t]}\Big(||v(\tau)||^2_{H^3}+||\vartheta(\tau)||^2_{H^3}\Big)&\leq C<\infty \quad\mbox{for all}\quad t\in(0,T^*),
\end{align*}
which implies that $T^*=+\infty$. This completes the proof of Theorem \ref{the1.1}. $\Box$

\section{Proof of Corollary \ref{cor1.1}}

In this section, we will give the proof of Corollary \ref{cor1.1}.

Case 1: d=2. Firstly, $W$ satisfys
\bbal
\pa_tW+\nu W=\pa_1\Theta =\pa_1\Theta_0e^{-\lambda t}, \quad W|_{t=0}=W_0=\nabla\times U_0.
\end{align*}
Formally,
\bbal
W=
\begin{cases}
e^{-\nu t}W_0+\frac{1}{\nu-\lambda}\pa_1\Theta_0(e^{-\lambda t}-e^{-\nu t}), \quad \nu\neq \lambda,\\
e^{-\nu t}W_0+te^{-\nu t}\pa_1\Theta_0,\quad \nu= \lambda.
\end{cases}
\end{align*}
Therefore, we can deduce that
\bbal
U=(-\Delta)^{-1}\na^{\bot}W=
\begin{cases}
e^{-\nu t}U_0+\frac{1}{\nu-\lambda}\pa_1(-\Delta)^{-1}\na^{\bot}\Theta_0(e^{-\lambda t}-e^{-\nu t}), \quad \nu\neq \lambda,\\
e^{-\nu t}U_0+te^{-\nu t}\pa_1(-\Delta)^{-1}\na^{\bot}\Theta_0,\quad \nu= \lambda£¬
\end{cases}
\end{align*}
where $\na^{\bot}=(\pa_2,-\pa_1)^T$.

\begin{lemma}\label{lem3.1} Let $d=2$. For small enough $\ep$, under the assumptions of Theorem \ref{the1.1}, the following estimates hold
\bal\label{estimate-fg}
E_0\leq C\ep||a_0||_{L^2}||\hat{a}_0||_{L^1}, \qquad F_0\leq C||\hat{a}_0||_{L^1}.
\end{align}
\end{lemma}
{\bf Proof of Lemma \ref{lem3.1}}\quad
With the above expressions of $U$ and $\Theta$, we can show that
\bbal
U\cd\na \Theta=
\bca
e^{-(\nu+\lambda) t}U_0\cd\na\Theta_0+\frac{1}{\nu-\lambda}(e^{-2\lambda t}-e^{-(\nu+\lambda) t})\pa_1(-\Delta)^{-1}\na^{\bot}\Theta_0\cd\na\Theta_0, \;\nu\neq \lambda\\
e^{-2\nu t}U_0\cd\na\Theta_0+te^{-2\nu t}\pa_1(-\Delta)^{-1}\na^{\bot}\Theta_0\cd\na \Theta_0, \;\nu=\lambda.
\eca\end{align*}
Notice that $U_0\cd\na\Theta_0=0$ and
\bbal
\pa_1(-\Delta)^{-1}\na^{\bot}\Theta_0\cd\na\Theta_0&=(\pa_2-\pa_1)\pa_1(-\Delta)^{-1}a_0\pa_1a_0
+\pa^2_1(-\Delta)^{-1}a_0(\pa_1-\pa_2)a_0.
\end{align*}
Using the classical Kato-Ponce product estimates and the fact the Fourier transform of a distribute belonging to $L^1$ lies in $L^\infty$, by a simple calculation, we obtain
\bbal
||\pa_1(-\Delta)^{-1}\na^{\bot}\Theta_0\cd\na\Theta_0||_{H^3}\leq C ||(\pa_2-\pa_1)a_0||_{L^\infty}||a_0||_{L^2}\leq C\ep||a_0||_{L^2}||\hat{a}_0||_{L^1},
\end{align*}
which implies
\bbal
||U\cd\na \Theta||_{H^3}\leq Ce^{-\lambda t}\ep||a_0||_{L^2}||\hat{a}_0||_{L^1}.
\end{align*}
Take the similar argument as the term $U\cd\na \Theta$, we also have
\bbal
||U\cd\na U||_{H^3}\leq Ce^{-\min\{\nu,\lambda\} t}\ep||a_0||_{L^2}||\hat{a}_0||_{L^1}.
\end{align*}
Thus, we complete the proof of Lemma \ref{lem3.1}. $\Box$\\

Case 2: d=3. Firstly, $W=\na\times U$ satisfys
\bbal
\pa_tW+\nu W=(\pa_2\Theta,-\pa_1\Theta,0)^\mathrm{T} =(\pa_2\Theta_0,-\pa_1\Theta_0,0)^\mathrm{T}e^{-\lambda t}, \quad W_{t=0}=W_0=\nabla\times U_0.
\end{align*}
Formally, we have
\bbal
W=
\begin{cases}
e^{-\nu t}W_0+\frac{1}{\nu-\lambda}(\pa_2\Theta_0,-\pa_1\Theta_0,0)^\mathrm{T} (e^{-\lambda t}-e^{-\nu t}), \quad \nu\neq \lambda,\\
e^{-\nu t}W_0+te^{-\nu t}(\pa_2\Theta_0,-\pa_1\Theta_0,0)^\mathrm{T} ,\quad \nu= \lambda.
\end{cases}
\end{align*}
Therefore, we can deduce that
\bbal
U&=(-\Delta)^{-1}\na\times W
\\&=
\begin{cases}
e^{-\nu t}U_0+\frac{1}{\nu-\lambda}(-\Delta)^{-1}
\big(\pa_1\pa_3\Theta_0,\pa_2\pa_3\Theta_0,-(\pa^2_1+\pa^2_2)\Theta_0\big)^\mathrm{T}(e^{-\lambda t}-e^{-\nu t}), \quad \nu\neq \lambda,\\
e^{-\nu t}U_0+te^{-\nu t}(-\Delta)^{-1}\big(\pa_1\pa_3\Theta_0,\pa_2\pa_3\Theta_0,-(\pa^2_1+\pa^2_2)\Theta_0\big)^\mathrm{T},\quad \nu= \lambda.
\end{cases}
\end{align*}

\begin{lemma}\label{lem3.2} Let $d=3$. For small enough $\ep$, under the assumptions of Theorem \ref{the1.1}, the following estimates hold
\bal\label{estimate-fg}
E_0\leq C\ep||a_0||_{L^2}||\hat{a}_0||_{L^1}, \qquad F_0\leq C||\hat{a}_0||_{L^1}.
\end{align}
\end{lemma}
{\bf Proof of Lemma \ref{lem3.2}}\quad
For the term $U\cd \na \Theta$, we can show that
when $\nu\neq \lambda,$
$$U\cd\na \Theta=
e^{-(\nu+\lambda) t}U_0\cd\na\Theta_0+\frac{1}{\nu-\lambda}(e^{-2\lambda t}-e^{-(\nu+\lambda) t})(-\Delta)^{-1}\big(\pa_1\pa_3\Theta_0,\pa_2\pa_3\Theta_0,-(\pa^2_1+\pa^2_2)\Theta_0\big)^\mathrm{T}\cd\na\Theta_0,
$$
while $\nu=\lambda$,
$$U\cd\na \Theta=e^{-2\nu t}U_0\cd\na\Theta_0+te^{-2\nu t}(-\Delta)^{-1}\big(\pa_1\pa_3\Theta_0,\pa_2\pa_3\Theta_0,-(\pa^2_1+\pa^2_2)\Theta_0\big)^\mathrm{T}\cd\na \Theta_0.$$
Notice that $U_0\cd\na\Theta_0=0$ and
\bbal
&\quad (-\Delta)^{-1}\big(\pa_1\pa_3\Theta_0,\pa_2\pa_3\Theta_0,-(\pa^2_1+\pa^2_2)\Theta_0\big)^\mathrm{T}\cd\na \Theta_0
\\&=(-\Delta)^{-1}\pa_1\pa_3\Theta_0\pa_1\Theta_0+(-\Delta)^{-1}\pa_2\pa_3\Theta_0\pa_2\Theta_0
-(-\Delta)^{-1}(\pa^2_1+\pa^2_2)\Theta_0\pa_3\Theta_0.
\end{align*}
Using the classical Kato-Ponce product estimates and the fact the Fourier transform of a distribute belonging to $L^1$ lies in $L^\infty$, after a simple calculation, we obtain
\bbal
&\quad||(-\Delta)^{-1}\big(\pa_1\pa_3\Theta_0,\pa_2\pa_3\Theta_0,-(\pa^2_1+\pa^2_2)\Theta_0\big)^\mathrm{T}\cd\na \Theta_0||_{H^3}
\\&\leq C ||\pa_3a_0||_{L^\infty}||a_0||_{L^2}\leq C\ep||a_0||_{L^2}||\hat{a}_0||_{L^1},
\end{align*}
which implies
\bbal
||U\cd\na \Theta||_{H^3}\leq Ce^{-\lambda t}\ep||a_0||_{L^2}||\hat{a}_0||_{L^1}.
\end{align*}
Take a similar argument as the term $U\cd\na \Theta$, we also have
\bbal
||U\cd\na U||_{H^3}\leq Ce^{-\min\{\nu,\lambda\} t}\ep||a_0||_{L^2}||\hat{a}_0||_{L^1}.
\end{align*}
Thus, we complete the proof of Lemma \ref{lem3.2}.

Now Corollary \ref{cor1.1} follows immediately from Lemma \ref{lem3.1} and Lemma \ref{lem3.2}. We
complete the proof of Corollary \ref{cor1.1} $\Box$\\

\section*{Acknowledgments} The work of Jinlu Li is supported by the National Natural Science Foundation of China (Grant No.11801090). The work of Weipeng Zhu is partially supported by the China National Natural Science Foundation under grant number 11901092 and Guangdong Natural Science Foundation under grant number 2017A030310634.

\end{document}